\documentclass[12pt,a4paper,leqno]{amsart}
\usepackage{latexsym,amssymb,amsthm,amscd}
\usepackage{eepic}
\usepackage{color}
\catcode`\@=11
\newcommand\lsection{\@startsection {section}{1}{\z@}%
                                   {-3.5ex \@plus -1ex \@minus -.2ex}%
                                   {1.0ex \@plus.2ex}%
                                   {\normalfont\large\bfseries}}
\catcode`\@=12
\newtheorem{thm}{Theorem}[section]
\newtheorem{cor}[thm]{Corollary}
\newtheorem{lem}[thm]{Lemma}

\theoremstyle{definition}

\newtheorem{rem}[thm]{Remark}

\newcommand{\NN}{\mathbf{N}}
\newcommand{\ZZ}{\mathbf{Z}}
\newcommand{\RR}{\mathbf{R}}

\newcommand{\BS}[1]{\boldsymbol{#1}}
\newcommand{\comp}{\text{\scriptsize$\circ$}}
\newcommand{\pd}[2]{\frac{\partial{#1}}{\partial{#2}}}

\newcommand{\ord}{\mathop{\mathrm{ord}}\nolimits}

\newcommand{\medsum}{\mathop{\textrm{\scriptsize $\mathop\sum$}}}

\newcommand{\medprod}{\mathop{\textrm{\small $\mathop\prod$}}}

\DeclareMathOperator{\im}{Im}

\newcommand{\AS}{\mathcal{AS}}
\newcommand{\inv}{^{-1}}
\date{March 2, 2010}

\title [Arc-analytic Homeomorphisms]{Inverse Function Theorems for Arc-analytic Homeomorphisms}

\thanks{Research partially supported by a Math\' ematiques en Pays de la Loire (MATPYL) grant.}
\author{Toshizumi Fukui, Krzysztof Kurdyka, Adam Parusi\'nski}

\address{Department of Mathematics, Faculty of Science, Saitama University, 255 Shimo- 
Okubo, Urawa 338, Japan}  

\email{tfukui@rimath.saitama-u.ac.jp}

\address{
Laboratoire de Math\'ematiques, UMR 5175 du CNRS, Universit\'e de Savoie,
Campus Scientifique, 73 376 Le Bourget--du--Lac Cedex, France,
}
\email{kurdyka@univ-savoie.fr}

\address { Laboratoire J. A. Dieudonn\'e 
U.M.R.  6621 du CNRS, Universit\'e de Nice Sophia-Antipolis,
28, Parc Valrose
06108 Nice Cedex 02, France} 

\email{adam.parusinski@unice.fr}

\begin{document}
\begin{abstract} 
We call a local homeomorphism $f: (\RR^n,0)\to(\RR^n,0)$ blow-analytic if it becomes real analytic, 
as a map, after composing with a finite number blowings-up with smooth nowhere dense centers.  
If the graph of $f$ is semi-algebraic then, by a theorem 
 of Bierstone and Milman, $f$ is blow-analytic if and only if 
it is arc-analytic: the image by $f$ of a parametrized 
real analytic arc $\gamma : (\RR , 0) \to  (\RR^n,0)$  is again a real analytic arc. 
  
We show  that if $f$ is blow-analytic, the inverse $f\inv$ of $f$ is Lipschitz, and the 
graph of $f$ is semialgebraic, then $f$ is Lipschitz and  $f \inv $ is blow-analytic.   The proof is by a 
motivic integration argument, using additive invariants on the spaces of arcs.  
\end{abstract}




\maketitle


\lsection{Introduction. }\label{Section:Introduction}

Let $M$, $N$ be real analytic manifolds. 
We say that $f:M\to N$ is {\it blow-analytic via $\pi$} if $\pi:\widetilde{M}\to M$
is a locally finite composition of
blowings-up with nonsingular nowhere dense centers and 
$f \comp \pi$ is analytic.
We say that $f$ is {\it blow-analytic} if there is such 
$\pi :\widetilde{M}\to M$,  that $f$ is blow-analytic via $\pi$.  We say that $f$ is semialgebraic if the graph 
of $f$ is semi-algebraic.

In this paper we show the following result.  

\begin{thm}\label{question}
Let  $f:(\RR^n,0)\to (\RR^n,0)$ be a semialgebraic homeomorphism such that $f$ is Lipschitz  and 
 $f\inv$ is blow-analytic.   
 Then  $f \inv$ is Lipschitz and $f$ is blow-analytic.    
\end{thm}

Theorem \ref{question}  gives a negative answer to Question 7.8 of  
 \cite{fukuipaunescu}.     As a corollary we obtain the following Inverse Function Theorem.  
 By  $C^\omega$ we mean real analytic.

 \begin{cor}\label{IFT}
Let  $f:(\RR^n,0)\to (\RR^n,0)$ be a semialgebraic homeomorphism such that 
 $f\inv$ is blow-analytic.    If $f$ is $C^k$, $k=1,2,\ldots, \infty, \omega$, then  so is $f \inv$.  
 \end{cor}
 
 The proof of Theorem \ref{question} uses the jet spaces of real analytic arcs 
and additive invariants of real algebraic sets.  
First we show the following theorem whose proof is given by 
a classical motivic integration argument on the jet spaces.

\begin{thm}\label{JacobianTh}
Let  $f:(\RR^n,0)\to (\RR^n,0)$ be a semialgebraic homeomorphism  such  that  
$f$ and $f^{-1}$ are blow-analytic.    
If the jacobian determinant $\det(df)$ is bounded, then there is a constant 
$c_1>0$ such that $c_1<|\det(df)|$. 
\end{thm}

\subsection{Blow-analytic and arc-analytic maps}

A map between real analytic manifolds  $f:M\to N$ is {\it arc-analytic}
 if for every real analytic 
arc germ $\gamma :(\RR ,0)\to M$ the composition $f\circ \gamma$ is analytic,
 see \cite{kurdykaAS}, \cite{aussoisKP}.  
A blow-analytic map is always arc-analytic. 
There is a partial reciprocal statement for $f$ subanalytic, 
see \cite{bm} and \cite{SubFun}, where the blow-analyticity is replaced by a similar notion 
 expressed in terms of local blowings-up.  
 In the semialgebraic case the blow-analiyicity and arc-analyticity are equivalent.

\begin{thm}\label{Arc-Blow} {\rm(Bierstone \& Milman,\cite{bm})}\\ 
If $M$ and $N$ are real algebraic manifolds and the graph of $f$ 
is semi-algebraic then the arc-analyticity of $f$ is equivalent to the blow-analyticity.  
Moreover, we may require that the blowings-up are along 
nonsingular real algebraic subvarieties. 
\end{thm} 

The blow-analytic maps were introduced in \cite{kuo}, 
see also surveys \cite{fukuikoikekuo}, \cite{fukuipaunescu}, in the context of blow-analytic equivalence 
of real analytic function germs, that is the equivalence induced by blow-analytic homeomorphisms.  
Neither blow-analytic equivalence implies the bi-lipschitz one, nor the vice-versa, cf. \cite{koikeparusinski2}, 
\cite{koikeparusinski3}.   Nevertheless  there is clear evidence that there is a relation between  
 blow-analylitic and Lipschitz 
property for  homeomorphisms, that should be further investigated and better understood.      
The following  Inverse Function Theorem for arc-analytic homeomorphism holds for 
homeomorphisms with subanalytic graphs, a class significantly bigger than the semialgebraic ones.

\begin{thm} \label{Thm:FKP} {\rm (\cite{FKP})}
Let  a subanalytic homeomorphism $f:(\RR^n,0)\to (\RR^n,0)$ be bi-Lipschitz and  arc-analytic. 
 Then $f\inv$ is also arc-anaytic. \end{thm}

Note : It is a widely accepted, \cite{fukuikoikekuo}, \cite{fukuipaunescu}, to call a homeomorphism $f:(\RR^n,0)\to (\RR^n,0)$,  \emph{a blow-analytic homeomorphism} if  both $f$ and $f\inv$ 
are blow-analytic,.  We avoid this terminology in this paper in order not 
to confuse it with a homeomorphism that is blow-analytic as a map (but maybe its inverse is not blow-analytic).

\subsection{Open problems.}

The proof of Theorem  \ref{question} and Theorem \ref{JacobianTh} is based on the motivic integration method  
on the space of real analytic arcs.  The essential point is the use of virtual Poincar\'e polynomial, 
an additive and multiplicative invariant of real algebraic varieties, that distinguishes their dimensions   
This part of the proof cannot be carried out, 
at the moment this paper is being written, to the subanalytic case, since such an invariant it is not known in this case.

We conjecture that Theorems \ref{question} and \ref{JacobianTh}  hold without 
assumption of semialgebraicity of the graph.  Then the graph has to be subanalytic since every 
blow-analytic map is subanalytic.  

We conjecture also that the following property holds :  \\
\emph{Let  $f:(\RR^n,0)\to (\RR^n,0)$ be a (semialgebraic) homeomorphism such that 
 $f\inv$ is blow-analytic.   If the jacobian determinant $\det(df)$ is bounded from above, then there is a constant 
$c_1>0$ such that $c_1<|\det(df)|$, and  $f$ is blow-analytic.    }


\lsection{Lipschitz and bi-Lipschiz maps. }\label{Section:Lipschitz}

In this section we show how Corollary  \ref{IFT} can be deduced from Theorem \ref{JacobianTh}.  
The argument is elementary.

Let $U$ be an open subset of $\RR^n$.  A map $f:U\to\RR^p$ is said to be {\it Lipschitz}
if there is a positive constant $L$ so that 
$$
|f(x)-f(y)|\le L |x-y|\qquad\forall x,y\in U.
$$
Let $U$ be an convex open subset of $\RR^n$ and let $f:U\to\RR$ be a continuous function with 
subanalytic graph.  Then there is an nowhere dense closed subanalytic subset $Z$ so that 
$f$ is analytic on $U-Z$.
\begin{lem}\label{Lem:Lip}
The function $f$ is Lipschitz if and only if all partial derivatives of $f$ are bounded on $U-Z$.
\end{lem}
\begin{proof}
If $f$ is Lipschitz, then the following inequality implies all 
directional derivatives are bounded whenever they exist. 
$$
\bigg|\frac{f(x+tv)-f(x)}{t}\Bigg|\le L|v|,\qquad v\in\RR^n. 
$$
Conversely, we assume that there is a positive constant $M$ so that 
$$
\biggl|\pd{f}{x_i}(x)\biggr|\le M, \qquad x\in U-Z, \ i=1,\cdots,n.
$$
For $x$, $x'\in U-Z$, we set $v=x'-x$ and write $v=(v_1,\dots,v_n)$. 
Then the mean value theorem implies that there is $\theta$ so that 
$$
f(x')-f(x)=\medsum_{i=1}^n\pd{f}{x_i}(x+\theta v)v_i, \quad0<\theta<1.  
$$
This makes sense when $x+\theta v\in U-Z$.  This implies that 
$$
|f(x')-f(x)|\le \medsum_{i=1}^n M |v_i| \le M \sqrt{n}\, |x'-x|. 
$$
Since $f$ is continuous, this means that $f$ is Lipschitz.
\end{proof}

\begin{cor}\label{JacBoundCor}
Let $f:(\RR^n,0)\to(\RR^n,0)$ be a subanalytic map. 
If $f$ is Lipschitz, then  $\det(df)$ is bounded. 
\end{cor}

\begin{cor}\label{BiLipCor}
Suppose that $f:(\RR^n,0)\to(\RR^n,0)$ is a subanalytic such that there are postive constants
 $c_1$, $c_2$ with
$$
c_1\le|\det(df)|\le c_2
$$
If $f$ is Lipschitz, then $f^{-1}$ is Lipschitz.  
\end{cor}

\begin{proof}
The jacobian matrix of $f^{-1}$ equals 
$$
\frac{1}{\det(df)}
(\textrm{cofactor matrix of the jacobian matrix of $f$})
$$
in the complement of a nowhere dense subanalytic set.  Since $f$ is Lipschitz, 
each coefficient of the jacobian matrix of $f^{-1}$ is bounded. 
We conclude by  Lemma \ref{Lem:Lip}.  
\end{proof}

Thus Corollary  \ref{IFT}  follows from Theorem \ref{JacobianTh} 
and Corollaries \ref{JacBoundCor}, \ref{BiLipCor}.



\lsection{Reduction to the normal crossing case.}  \label{reduction}

A real modification  is a classical notion introduced in  \cite{kuo}.  It is a natural generalization of 
the notion of  modification in algebraic geometry.  For instance,   a locally finite composition of
blowing-ups with nonsingular nowhere dense centers is a real modification 
in the sense of  \cite{kuo}.  In was shown in \cite{koikeparusinski2}, \cite{fukuipaunescu},  
 that every real modification satisfies {the unique lifting of generic arc property: 
 
\begin{thm} {\rm  (\cite{koikeparusinski2}, \cite{fukuipaunescu})}
 Let $\tau : M\to N$  be a real modification.  Then 
there is a closed subanalytic nowhere dense 
$A\subset N$ such that for every real analytic curve germ 
 $\gamma : (\RR ,0) \to (N,p)$, if the image of $\gamma$ is not entriely contained in $A$ then 
 there is a unique  real analytic  lift  $\tilde \gamma : (\RR ,0)  \to (M,\tilde p)$,  
such that  $\tau \circ \tilde \gamma = \gamma$.  
\end{thm}

Suppose that, as in the assumption of Theorem \ref{JacobianTh}, the mapping 
 $f:(\RR^n,0)\to(\RR^n,0)$  and its inverse $f^{-1}$ are blow-analytic.   
 Then by argument of \cite{kuo}, proof of Proposition 2, 
there is a real analytic manifold $M$, and real modifications
$\sigma: (M,\sigma \inv (0) )\to (\RR^n,0)$, $\sigma' : (M, {\sigma'}\inv (0))\to (\RR^n,0)$,  
such that $f\circ \sigma =\sigma'$.  

Let $\pi$ be a sequence of blowing-ups with non-singular centers, so that
$\det(d\sigma )\circ \pi$, $\det(d\sigma' )\circ \pi$ 
are simultaneously normal crossings.
We may 
assume that the new centres are in normal crossings with all old exceptional
divisors.  Then  the jacobian of $\pi$,
$\det(d\pi)$, is normal crossing.  Therefore  the jacobian determinants  
$\det d(\sigma \circ \pi) = ( (\det(d\sigma))\circ \pi)\cdot   \det(d\pi)$, 
$\det d(\sigma' \circ \pi) = ( (\det(d\sigma'))\circ \pi)\cdot   \det(d\pi)$, 
are normal crossings and, moreover,  after composing with blowing-ups
with nonsingular centers, if necessary,
we may assume that they are normal crossings 
in the same system of coordinates.

We thus may assume that the critical loci of $\sigma$ and $\sigma'$
are simultanuously normal crosssing and denote them respectively  by 
\begin{equation}\label{multiplicities}
 E=\sum_{i\in I}\nu_i E_i,\quad \textrm{ and }\quad 
 E' =\sum_{i\in I'}\nu'_i E_i .
\end{equation}
Since 
$(\det(df)\circ \sigma)\cdot\det(d\sigma)= \det(d\sigma')$, 
 for a generic real analytic arc $\gamma$ at a generic 
$p\in E_i$
$$
  \ord_0 (\det(df)\circ \sigma\circ \gamma)=\nu'_i-\nu_i.
$$
Thus if   the jacobian determinant $\det(df)$ is bounded, then 
$\nu'_i\ge\nu_i$ for all $i\in I$.   
Under assumptions of Theorem \ref{JacobianTh}, we show in the next section that 
this implies that  $\nu'_i=\nu_i$ for all $i\in I$.



\lsection{Arc spaces and additive invariants.}\label{MainSection}

\subsection{Constructible sets}
For the proof we need \emph{the virtual Poincar\'e polynomial}, 
an additive and multiplicative invariant, 
first introduced for real algebraic varieties in \cite{virtual}, 
and then for a wider class of  $\AS$ sets in \cite{fichou} 
and \cite{weight}.   A semialgebraic 
subset $X$ of a compact real algebraic variety $V$ is called \emph{an $\AS$ set} 
if it is a finite set-theoretic combination of semi-algebraic arc-symmetric subsets of  $V$, 
cf. \cite{kurdykaAS},  \cite{aussoisKP}.    \emph{The virtual Poincar\'e polynomial of $X$} 
$$\beta (X) = \sum \beta_i(X) u^i \in \ZZ [u]$$
 satisfies the following properties, see \cite{virtual}, \cite{fichou}, 
\begin{enumerate}
\item
\emph{Additivity:} For finite disjoint union $X= \sqcup X_i$, $\beta(X) = \sum \beta (X_i)$. 
\item
\emph{Multiplicativity:}  $\beta (X\times Y) =  \beta(X)  \beta (Y)$.  
\item
\emph{Degree:}  For $X\ne \emptyset$, $\deg \beta (X) =  \dim X$ 
and the leading coefficient $\beta(X)$ is  strictly positive.  
\end{enumerate}
If $X$ is compact and nonsingular then $\beta_i(X) = \dim H_i(X;\ZZ_2)$. 

 The virtual Poincar\'e polynomial is an invariant of real analytic isomorphisms with 
 semialgebraic graph, see  \cite{fichou}, 
 and, more generally,  of bijections with $\AS$ graphs, see \cite{weight}.  
 
 In the rest of this section we mean by \emph{constructible set}, an $\AS$ set, 
and by \emph{constructible map},   a map with constructible graph.  
 We refer the reader to   \cite{aussoisKP}, and also to \cite{weight}, \cite{mathannalen}, 
 for more precise discussion. 
 
 By \emph {morphism} we mean a Nash map $f:M\to N$ that is a real analytic map with 
semi-algebraic graph.  For our 
 purpose we may suppose that $M$ and $N$ are nonsingular real algebraic varieties, 
but, in general, we  may consider $M$ and $N$ to be Nash manifolds, \cite{BCR}.  
 By a \emph{modification} we mean a real modification that is a Nash map, 
i. e. a Nash modification in the sense of 
 \cite{fichou2}.  A regular proper birational map is a standard example of such modification.


\subsection{Arc spaces}\label{arcspaces}
We use the technique developed in \cite{DL}, and adapted to the real analytic set-up in 
\cite{KoikeParusinski}, \cite{fichou}, \cite {fichou2}.    
Let $M$ be a real analytic manifold and let $S$ be  a subset of $M$.
Consider the  arc space
$$
\mathcal L(M,S):=\{
\gamma : (\RR,0) \to (M,S) , \
\textrm{analytic}\}.
$$
For a real analytic map $\sigma:M\to\RR^n$, set
$$
\mathcal L:={\mathcal L}(\RR^n,0),\quad
\widetilde{\mathcal L}:=\mathcal L(M,\sigma^{-1}(0)),\quad
\widetilde{L}_k:=\mathcal \bigcup_{x\in\sigma^{-1}(0)} L_k(M,x),
$$
where $L_k(M,x)$ denotes the set of $k$-jets of elements of $\mathcal L(M,x)$.
Setting $L_k=L_k(\RR^n,0)$,  we have the following commutative diagram of natural maps:
$$
\CD
\widetilde{\mathcal L}@>{\sigma_*}>
> \mathcal L\\
@V{p_k}VV@VV{p_k}V\\
\widetilde{L}_k@>{\sigma_{*,k}}>>L_k
\endCD$$
where $p_k$ denote the maps defined by taking the $k$-jets.  Consider 
$$
\mathcal B_e(\sigma)=\{\gamma\in\widetilde{\mathcal L}
: \ord_{\gamma}\det(d\sigma)=e\},   \qquad 
B_{k,e}(\sigma)=p_k(\mathcal B_e(\sigma)), 
$$
where $ \ord_{\gamma}\det(d\sigma)$ is defined as the order of $ \det(d\sigma) (\gamma (t))$ at 
$t=0$.  
It is clear that $B_{k,e}(\sigma)$ is a difference of two analytic sets.  
By Lemma 2.11 of \cite{fichou2} we have the following.

\begin{lem}\label{Lem:Tri}
Let $\sigma:(M,E_0)\to(\RR^n,0)$ be a Nash modification.  
Assume  $k\ge2e$. Then $\sigma_{*,k}(B_{k,e}(\sigma))$ is constructible and
 $B_{k,e}(\sigma)\to \sigma_{*,k}(B_{k,e}(\sigma))$
is a piecewise trivial fibration  with fiber $\RR^e$. 
\end{lem}

\begin{rem}
\item[ (i)]
The statement of Lemma \ref{Lem:Tri} means that there is a finite partition of 
$\sigma_{*,k} (B_{k,e}(\sigma))$ into constructible  sets so that  
$B_{k,e}(\sigma)\to \sigma_{*,k}(B_{k,e}(\sigma))$ over each piece is isomorphic  
to a trivial fibration by a constructible homeomorphism.  
\item[ (ii)]
Lemma \ref{Lem:Tri} does not hold, in general, if one assumes only that $\sigma$ is  a proper 
Nash map.  In the statement of  Lemma 2.11 of \cite{fichou2} the assumption   
that $h$ is a Nash modification is missing.
\end{rem}


\smallskip
Let $\sigma:(M,E_0)\to(\RR^n,0)$, $E_0 =\sigma^{-1}(0)$,  be a Nash modification.
Assume that the critical locus of $\sigma$,  $E := \{ \det(d\sigma)=0\}$,   
 is a divisor with normal crossings.  
Let 
$$
(\det(d\sigma))_0=\medsum_{i\in {I}}\nu_i E_i  , 
$$
where $E_i$ are components of $E$, and $\nu_i >0$ for each $i\in  I$.  Denote 
$\nu_{\max}= \max\{\nu_i:i\in {I}\}$.   
We also assume that $E_0=\sigma \inv (0)$ is a union of components of $E$ 
$$
\sigma \inv (0) =\bigcup _{i\in {I_0}}  E_i  .
$$
For a vector $\BS{j}=(j_i)_{i\in {I}}$, $j_i\in \NN$, 
we set $J =J(\BS {j} )=\{i:j_i\ne0\}\subset I$, 
$E_J = \bigcap_{i\in J} E_i$ and  
$\mathring E_J=\bigcap_{i\in J}E_i\setminus\bigcup_{j\in I\setminus J}E_j$.  
We only consider such $\BS {j}$  that $\sigma (E_J) =\{0\}$.  
For such $\BS {j}$  we denote
$$
\mathcal B_{\BS{j}}
=\{\gamma\in\mathcal L
(M, \mathring E_I):\ord_\gamma E_i=j_i, \ i\in J\} \subset \widetilde {\mathcal L} 
$$
and for $k \in \NN$ 
$$
B_{k,\BS{j}}:=p_k(\mathcal B_{\BS{j}}), \qquad X_{k,\BS{j}}(\sigma)= \sigma_{*,k}(B_{k,\BS{j}}) .
$$ 
Finally we set 
\begin{align*}
A_k(\sigma)=& \{\BS{j}  : \ \sigma(E_J)=\{0\} \text { and } \langle\BS{\nu},\BS{j}\rangle\le k/2\},
\end{align*}
where $\langle \BS{\nu},\BS{j}\rangle :=\sum_{i\in {I}}\nu_ij_i $.

\medskip
\begin{lem}\label{Partition}
The sets $X_{k,\BS{j}}(\sigma)$, $\BS{j}\in A_k(\sigma)$,  
are constructible subsets of 
$L_k$ and $\dim X_{k,\BS{j}}(\sigma)=
n(k+1)-s_{\BS{j}}-\langle\BS{\nu},\BS{j}\rangle$,   
where $s_{\BS{j}}=\sum_{i\in J} j_i$.  
We have a disjoint union 
$$
L_k=Z_k(\sigma)\sqcup
\bigsqcup_{\BS{j}\in A_k(\sigma)}X_{k,\BS{j}}(\sigma) ,
$$
and the constructible set  $Z_k(\sigma)$ satisfies 
$\dim Z_k(\sigma)<n(k+1)-k/2\nu_{\max}$.
\end{lem}

\begin{proof}
Fix $\BS{j}$ such that $\mathring E_J\ne \emptyset$, $\sigma(E_J) = \{0\}$, and 
$0\le j_i\le k+1$ for $i\in I$.   Since the fiber of the natural projection 
$B_{k,\BS{j}} ( \mathring E_J)\to \mathring E_J$ is
$$
\medprod_{i\in J}(\RR^*\times\RR^{k-j_i})\times(\RR^k)^{n-|J|}
\simeq
(\RR^*)^{|J|}\times\RR^{nk-s_{\BS{j}}},
$$
we conclude that 
$$
\dim B_{k,\BS{j}}= n(k+1)-s_{\BS{j}}.
$$

Assume that $\BS{j} \in A_k(\sigma)$.  The sets $X_{k,\BS{j}}(\sigma)$ are constructible.  
 Indeed, 
$X_{k,\BS{j}}(\sigma)$ is the image of a constructible set $B_{k,\BS{j}}$ and, by Lemma 
\ref{Lem:Tri},   $B_{k,\BS{j}}  \to X_{k,\BS{j}}(\sigma)$ has all the fibers of constant 
Euler characteristic with compact supports equal to $\pm 1$.  Therefore,  the characteristic function 
$\mathbf 1 _{X_{k,\BS{j}}(\sigma)} $ is Nash constructible in the sense of \cite{aussoisKP}, 
see also \cite {ens} Section 5 or \cite {weight}, 
that implies that   $X_{k,\BS{j}}(\sigma)$ is constructible.  Its dimension is given by 
$$
\dim X_{k,\BS{j}}  (\sigma)= n (k+1)-s_{\BS{j}}-e.  
$$
We also have 
$$
X_{k,\BS{j}}(\sigma)\cap X_{k,\BS{j}'}(\sigma)
=\emptyset\qquad\textrm{if}\qquad \BS{j}\ne\BS{j}', \quad
\BS{j},\BS{j}'\in A_k(\sigma).  
$$

Assume that $\BS{j}\not \in A_k(\sigma)$, that is $k<2\langle \BS{\nu},\BS{j}\rangle $. Then 
$$
\dim X_{k,\BS{j}} (\sigma) \le \dim B_{k,\BS{j}}
(k+1)-s_{\BS{j}}.
$$

Since
$k/2< \langle\BS{\nu},\BS{j}\rangle \le\nu_{\max}s_{\BS{j}}$,
we have 
$$
\dim X_{k,\BS{j}} (\sigma)
<n(k+1)-\frac{k}{2\nu_{\max}}, 
$$
as claimed.  
\end{proof}

\begin{cor}\label{Cor:VBN}
By Lemma \ref{Lem:Tri},  
$\beta(X_{k,\BS{j}}(\sigma))=
\beta( \mathring E_{J})(u-1)^{|J|}
u^{nk-s_{\BS{j}}-\langle\BS{\nu},\BS{j}\rangle}$ and hence 
$$u^{nk}=
\beta(Z_k(\sigma))+\sum_{\BS{j}\in A_k(\sigma)}
\beta( \mathring E_{J})(u-1)^{|J|}
u^{nk-s_{\BS{j}}-\langle\BS{\nu},\BS{j}\rangle}. 
$$
\end{cor}

\medskip
\begin{thm}\label{Thm1}
Let $\sigma:(M,E_0)\to(\RR^n,0)$ and $\sigma':(M',E' _0)\to(\RR^n,0)$ 
are two Nash modifications.  
Suppose that the critical loci $E$,  $E'$, of $\sigma$ and  $\sigma'$ 
are normal crossing divisors 
$$
(\det(d\sigma))_0=\medsum_{i\in {I}}\nu_iE_i,\qquad
(\det(d\sigma'))_0=\medsum_{i\in {I'}}\nu'_iE'_i , 
$$
and  that $\sigma \inv (0)$, 
respectively $(\sigma') \inv (0)$, is a union of components of $E$, resp. of $E'$.  

Let $F:M\to M'$ be a Nash isomorphism such that $F(\sigma \inv (0)) = {\sigma'}\inv (0)$ 
and $F(E_i) = E'_{\varphi(i)}$, for  $i\in {I}$ and $\varphi: {I}\to {I}'$ is a bijection.  
If $\nu_i\le\nu'_{\varphi(i)}$ for all $i\in {I}$, 
then $\nu_i=\nu'_{\varphi(i)}$ for all $i\in {I}$. 
\end{thm}

\begin{proof}
We identify $M'$ with $M$, $E'_{\varphi(i)}$ with $E_i$, and ${I}'$ with ${I}$, 
and denote them by the same letters.   
By Corollary \ref{Cor:VBN}  
\begin{equation}\label{B}
P=Q'-Q+\beta(Z_k(\sigma'))-\beta(Z_k(\sigma))
\end{equation}
where 
\begin{align*}
P = &\sum_{\BS{j}\in A_k(\sigma)\cap A_k(\sigma')}
\beta (\mathring E_J)(u-1)^{|J|}
u^{nk-s_{\BS{j}}-\langle\BS{\nu}',\BS{j}\rangle}
(u^{\langle\BS{\nu}'-\BS{\nu},\BS{j}\rangle}-1),  \\
Q =&\sum_{\BS{j}\in A_k(\sigma)\setminus A_k(\sigma')}
\beta  ( \mathring E_J)(u-1)^{|J|}
u^{nk-s_{\BS{j}}-\langle\BS{\nu},\BS{j}\rangle}, \\  
Q' = &\sum_{\BS{j}\in A_k(\sigma')\setminus A_k(\sigma)}
\beta ( \mathring E_J)(u-1)^{|J|}
u^{nk-s_{\BS{j}}-\langle\BS{\nu}',\BS{j}\rangle} .
\end{align*}
The assumption $\nu_i\le \nu'_{i}$ gives $A_k(\sigma) \supset A_k(\sigma ' )$ and therefore 
$Q' \equiv 0$.

Let
$$
C_k=\{s_{\BS{j}}+\langle\BS{\nu},\BS{j}\rangle:
\BS{j}\in A_k(\sigma'),\ 
\langle\BS{\nu}'-\BS{\nu},\BS{j}\rangle>0\} 
$$
and suppose that for  $k$ big enough $C_k$ is nonempty.  The minimum $c_k= \min C_k$ stabilizes.  
Thus denote $c=c_k$ for $k$ big enough, say $k\ge k_0$.   Then, for $k\ge k_0$,  
$$
\deg P = \max \{ \dim B_{k, \BS {j}}; \BS {j} \in  A_k(\sigma') \} = 
n(k+1)-c.  
$$

But this, for $k$ big enough, contradicts the following  lemma. 

\begin{lem}\label{Lem:bounds}
We have the following degree bounds :
\begin{align*}
\deg \beta(Z_k(\sigma))<&n(k+1)-\frac{k}{2\nu_{\max}},&
\deg \beta(Z_k(\sigma'))<&n(k+1)-\frac{k}{2\nu'_{\max}}, \\
\deg Q<&n(k+1)-\frac{k}{2\nu'_{\max}},&
\hfil  
\end{align*}
where 
$\nu_{\max}=\max\{\nu_i:i\in  {I}\}$, 
$\nu'_{\max}=\max\{\nu'_i:i\in I'\}$.
\end{lem}

Indeed, the degree bounds for 
$\beta(Z_k(\sigma))$, $\beta(Z_k(\sigma'))$ are a 
consequence of Lemma \ref{Partition}. 
If $\BS{j}\in A_k(\sigma)\setminus A_k(\sigma')$, then 
 $k/2<\langle\BS{\nu}',\BS{j}\rangle\le\nu'_{\max}s_{\BS{j}}$. 
We thus have 
$$
s_{\BS{j}}+\langle\BS{\nu},\BS{j}\rangle
\ge s_{\BS{j}}
> \frac{k}{2\nu'_{\max}}
$$
This implies the degree bound for $Q$ and  ends the proof of Theorem \ref{Thm1}.  
\end{proof}


\subsection{Proof of Theorem \ref{JacobianTh}}
If $f$ is arc-analytic and its graph is semialgebraic then, by \cite{bm}, $f$
 is blow-analytic via a sequence of blowings-up  
with nonsingular algebraic centers.  Similarly for $f\inv$.  Therefore, by \cite{kuo}, 
 there are Nash modifications $\sigma, \sigma' : M\to \RR^n$ such that $f\circ \sigma =\sigma'$.  
Consequently, Theorem   \ref{JacobianTh} follows from Theorem \ref{Thm1} and section \ref{reduction}.  \qed

\begin{rem}\label{fukui}
Theorem \ref{Thm1} holds true if we require $F : M\to M'$ to be only a homeomorphism satisfying
 $F(\sigma \inv (0)) = {\sigma'}\inv (0)$ and $F(E_i) = E'_{\varphi(i)}$.  Indeed, 
in this case $\beta (E_i) = \beta (E'_{\varphi(i)})$ and hence by 
additivitty $\beta  ( \mathring E_J )= \beta  ( \mathring E_{\varphi (J)})$.  

Suppose that $f$ is given by a commutative diagram  
$$
\CD
(M,E_0)@>{F}>
> (M',E'_0)\\
@V{\sigma }VV@VV{\sigma ' }V\\
(\RR^n,0) @>{ f}>> (\RR^n,0)
\endCD$$
with $\sigma $ and $\sigma '$ Nash modifications  (e.g. compositions of regular algebraic
 blowings-up) and $F$ a homeomorphism satisfying the above mentioned properties. 
In this case, if the jacobian determinant $\det(df)$ is bounded, then there is a constant 
$c_1>0$ such that $c_1<|\det(df)|$. 
In particular,  if  $f$ is Lipschitz then so is its inverse  $f^{-1}$. 

This shows that Theorems \ref{Thm1} and  \ref{JacobianTh}
 hold without the assumption of semialgebraicity in the two variable case.  
\end{rem}


\lsection{Proof of Theorem \ref{question}}

We may suppose that there is a composition of blowings-up with non-singular nowhere 
dense centers $\sigma : (M,E_0)\to (\RR^n,0)$
  such that  $\sigma' = f \circ \sigma  :(M,E_0)\to (\RR^n,0)$ is analytic.  
 Moreover, by performing additional blowings-up, as in section 
\ref{reduction}, we may assume that 
 the critical loci of $\sigma$ and $\sigma '$ are normal crossings as in 
\eqref {multiplicities}.   We also assume that that $E_0 = \sigma\inv (0) = {\sigma'}\inv (0)$
 is a union of components of $E$.

We denote by $ \mathcal L (\sigma)$ and  $ \mathcal L (\sigma')$ the space of 
real analytic arcs of the source and, respectively, of  the target of $f$.  
Thus both are equal to  ${\mathcal L}(\RR^n,0)$.  Since $f$ is arc-analytic it induces a map 
$$f_{*} : \mathcal L (\sigma) \to \mathcal L (\sigma' ). $$
 In general, this map does not factor to a map  $f_{*,k} : L_k (\sigma) \to L_k (\sigma' ) $.  This is the case, for 
 all $k$, if  $f_{*}$ preserves the order of contact of parametrized curves.  Therefore, 
 by the curve selection lemma,  $f_{*,k}$ is well defined for all $k\in \NN$ if and only if $f$ is Lipschitz,  
 that we would like show.  

Note that $f_{*}$ is injective since $f$ is a homeomorphism, and therefore 
$f\inv _{*}$ is well-defined on the image of $\sigma'_*$.  
Moreover, $f\inv$ is Lipschitz  and therefore the jacobian determinant of $f\inv $ is bounded and 
$f\inv _{*}$ factors to the jet spaces.  
More precisely we have the following 
result.

\begin{lem}\label{Lem:LipConcl}
We have $\nu'_i \le \nu_i$ for all $i\in I$. 
Moreover, let $\gamma_1, \gamma_2 \in \widetilde L^k$.  
 If  $ \sigma'_{*,k} (\gamma_1 ) =  \sigma'_{*,k} (\gamma_2 ) $ 
then 
$ \sigma_{*,k} (\gamma_1 ) =  \sigma_{*,k} (\gamma_2 ) $.  
Hence, $f\inv_{*,k}$ is well-defined on $\im (\sigma'_{*,k}) $.  
\end{lem}

Consider the induced maps  on $k$-jets  
$$
\CD
\widetilde{L}_k  @= \widetilde{L}_k \\ 
@V{\sigma_{*,k}}VV@VV{{\sigma'_{*,k}}}V \\
{L}_k (\sigma ) @<{f\inv_{*,k}}<<  \im (\sigma'_{*,k})
\endCD$$

Then, similarly to Subsection \ref{arcspaces} we have a decomposition 
$$
\im (\sigma'_{*,k})  = \ Z_k(\sigma' )\sqcup
\bigsqcup_{\BS{j}\in A_k(\sigma)}   X_{k,\BS{j}}(\sigma').  
$$
Note that  $X_{k,\BS{j}}(\sigma ' ) = ( f \inv _{*,k}) \inv  (X_{k,\BS{j}}(\sigma))  
= \sigma'_{*,k} (B_{k,\BS{j}} )$, and 
$Z_k (\sigma ' )=  (f\inv  _{*,k}) \inv  (Z _k (\sigma))$.  
We have $\sigma_{*,k} =  f \inv _{*,k} \circ \sigma'_{*,k} $ and since  $\sigma_{*,k} $ is epi so is  $  f \inv _{*,k} $.  
Therefore 
$$ \dim  Z_k (\sigma')\ge \dim Z _k (\sigma) .$$ 

By Lemma \ref{Lem:LipConcl}, $A_k (\sigma') \supset A_k (\sigma)$.   
We divide $A_k(\sigma)$ into two pieces  
\begin{align*}
A'_k(\sigma)=\{\BS{j} \in A_k(\sigma)   :  \langle\BS{\nu},\BS{j} \rangle =  \langle\BS{\nu'},\BS{j}  \rangle \}, \\
A''_k(\sigma)= \{\BS{j} \in A_k(\sigma)   :  \langle\BS{\nu},\BS{j} \rangle >  \langle\BS{\nu'},\BS{j} \rangle \}.
 \end{align*}
If $\sigma'$ is not a Nash modification then the statement of Lemma \ref{Lem:Tri} may not hold for $\sigma'$.  Nevertheless, the computations of the proofs of Lemma 4.2 \cite{KoikeParusinski} or Lemma 3.4 of \cite {DL}  give
 the following local result.  We denote by $p_{k,k-e} : \tilde  L_k \to \tilde L_{k-e}$ the truncation. 

\begin{lem}\label{Lem:LocTri:new}
Assume  $k\ge2e$.  Let $\gamma \in B_{k,e}(\sigma') $,  $\gamma(0)  \in E_0$.   
Then the fibre of 
$$
 p_{k,k-e}\inv  (p_{k,k-e}(\gamma)) \cap B_{k,e}(\sigma')\to \sigma'_{*,k} (B_{k,e}(\sigma')),   
$$
containing $\gamma$,  is an affine subspace of $p_{k,k-e}\inv  (p_{k,k-e}(\gamma))$ of dimension $e$.  
\end{lem}

This lemma holds without assuming that $\sigma '$ is birational or a Nash modification.  Thus any fibre of 
$B_{k,e}(\sigma') \to \sigma'_{*,k} (B_{k,e}(\sigma'))$ contains an affine space of dimension $e$.  

\begin{cor}
If  $\BS{j} \in A_k(\sigma)$ then 
$\dim   X_{k,\BS{j}}(\sigma' ) =  \dim X_{k,\BS{j}}(\sigma)  + \langle \BS{\nu} - \BS{\nu'},\BS{j} \rangle$. 
Moreover, if  $\BS{j} \in A'_k(\sigma)$, then   $f\inv_{*,k}$ induces a bijection 
$X_{k,\BS{j}}(\sigma ) \leftarrow  X_{k,\BS{j}}(\sigma' )$.  
\end{cor} 

\begin{proof}
Let $\gamma \in B_{k,\BS{j}}$.  Then we have the inclusions
\begin{align*}
{\sigma'_{*,k}} \inv (\sigma'_{*,k}(\gamma )) \subset 
\sigma_{*,k} \inv ( \sigma_{*,k}(\gamma )) \simeq \RR^e \subset  p_{k,k-e}\inv  (p_{k,k-e}(\gamma)),
\end{align*}
where $e = \langle\BS{\nu},\BS{j} \rangle$.  Therefore,  by Lemma \ref{Lem:LocTri:new},  
$$
{\sigma'_{*,k}} \inv (\sigma'_{*,k}(\gamma ))  \simeq \RR^{e'} ,
$$ 
where $e' = \langle\BS{\nu'},\BS{j} \rangle$.  This gives the first claim of corollary.  

If  $\BS{j} \in A'_k(\sigma)$, then $ {\sigma'_{*,k}} \inv (\sigma'_{*,k}(\gamma )) = 
\sigma_{*,k} \inv ( \sigma_{*,k}(\gamma ))$ and the second claim of corollary follows.   
\end{proof}

The map $f\inv_{*,k}$ is constructible, that is its graph is an $\AS$ set, and therefore, by  \cite {weight},  
$$\beta (X_{k,\BS{j}}(\sigma ) ) = \beta (\tilde X_{k,\BS{j}}(\sigma' )) .$$
Thus we have 
\begin{align}\label{beta}
&  \beta (L_k) - \beta(\im (\sigma'_{*,k} ) ) \hfill \\ \nonumber 
& = \sum_{\BS{j} \in A''_k (\sigma)}
(\beta (X_{k, \BS{j}} (\sigma) ) -  \beta ( \tilde X_{k, \BS{j}} (\sigma ')) ) + (\beta (Z_k(\sigma)) -\beta(
\tilde Z_k(\sigma ' )).
\end{align}
The leading coefficient of the left-hand side  is positive (if $\im (\sigma'_{*,k} )\ne L_k$), and the leading 
coefficient of the first summand of the right-hand side is negative.  The leading 
coefficient of the second summand is also negative unless  $\dim \tilde Z_k (\sigma')= \dim Z _k (\sigma)$.  

Thus, necessarily,   $\dim \tilde Z_k (\sigma')= \dim Z _k (\sigma)$ and  if  $\BS{j}\in A''_k(\sigma)$ then 
 $\dim \tilde X_{k,\BS{j}}(\sigma ') \le  \dim Z _k (\sigma)$.  But this is impossible by the argument 
 of  proof  of Theorem  \ref{JacobianTh}.  Indeed, for fixed  $\BS{j}$ by letting $k\to \infty$, we obtain by lemma 
 \ref {Lem:bounds} the opposite inequality  $\dim \tilde X_{k,\BS{j}}(\sigma ') >  \dim Z _k (\sigma)$.  Thus 
 $ \nu_i = \nu'_i   \text { for all } i\in I$ and $f$ is Lipschitz by Corollary \ref{BiLipCor} .  
 
 Finally, $f\inv$ is blow-analytic by Theorem \ref{Thm:FKP}.  
\qed

\begin{rem}
In the proof of Theorem \ref{question} it is not necessary to use  Theorem \ref{Thm:FKP}.  Instead, 
we can argue directly as follows.  By \eqref{beta} we see that $\sigma'_{*,k}$ is surjective,
 since the codimension of 
$Z_k(\sigma')$ goes to infinity as $k\to \infty$.   If $f\inv $ were not blow-analytic, 
there would have been a real analytic arc  $\gamma(t) = \sum a_i t^i$ such that 
\begin{align}\label{notanalytic}
f\inv (\gamma (t))= \sum _{i=1}^m b_i t^i + b t^{p/q} + \cdots , \quad b\ne 0
\end{align}
where $m<p/q<m+1$.  Changing the higher terms of $\gamma$, if necessary, we may assume that for  $k\gg m$, 
$p_k (\gamma) \in X_{k,\BS{j}}(\sigma') $, and $f\inv_{*,k}(\gamma) \in X_{k,\BS{j}}(\sigma )$, 
that contradicts \eqref{notanalytic}.  

\end{rem}


\bigskip
\begin{thebibliography}{99}


\bibitem{bm}  E. Bierstone, P. D. Milman,
\emph{Arc-analytic functions,}
Invent. math.  \textbf{101} (1990), 411--424.

\bibitem {BCR}  J. Bochnak, M. Coste, M.-F. Roy,
\emph{Real Algebraic Geometry}, Springer Verlag, New York,
1992.

\bibitem{DL}  J. Denef, F. Loeser,
\emph{Germs of arcs on singular algebraic varieties and motivic integration},
  Invent. math. \textbf{135} (1999) 201--232.


\bibitem{fichou}  G. Fichou,
\emph{Motivic invariants of arc-symmetric sets and blow-Nash equivalence},
  Compositio Math. \textbf{141} (2005) 655--688.

\bibitem{fichou2}  G. Fichou,
\emph{ Zeta functions and Blow-Nash equivalence,}
  Annales Polonici Math., special volume in memory of  S. 
{\L}ojasiewicz, \textbf{87} (2005) 111--126

\bibitem{fukuikoikekuo} T. Fukui, S. Koike, T.-C. Kuo, 
  \emph{Blow-analytic equisingularities, properties, problems
  and progress}, Real Analytic and Algebraic Singularities
  (T. Fukuda, T. Fukui, S. Izumiya and S. Koike, ed),
  Pitman Research Notes in Mathematics Series,
  \textbf{381} (1998), pp. 8--29.
  
  \bibitem{FKP} 
T.~Fukui, K.~Kurdyka,  L.~Paunescu,
\emph{An inverse mapping theorem for arc-analytic homeomorphisms, } 
Geometric singularity theory, 49--56, Banach Center Publ., 65, Polish Acad. Sci., Warsaw, 2004.


\bibitem{fukuipaunescu} T. Fukui, L. Paunescu, 
  \emph{On blow-analytic equivalence}, in 
"Arc Spaces and Additive Invariants in Real Algebraic Geometry",  
Panoramas et Synth\`eses \textbf{24}, 2007, S.M.F., pp. 87--122


\bibitem{kkf} 
T.~Fukui, S.~Koike,   T.-C.~Kuo,
Blow-analytic equisingularities, properties, problems and progress, 
in ``Real analytic and algebraic singularities'',
Pitman Research Notes in Mathematics Series,
{\bf 381}, 1997, Longman, 8--29.


\bibitem{KoikeParusinski}
S.~Koike and A.~Parusi\'nski,  
 \emph{Motivic-type invariants of blow-analytic equivalence},
  Ann. Inst. Fourier \textbf{53} (2003), 2061--2104.


\bibitem{koikeparusinski2} S. Koike, A. Parusi\'nski,
 \emph{Blow-analytic equivalence of two variable real analytic function germs}, 
arXiv:0710.1046, to appear in J. Algebraic Geometry

\bibitem{koikeparusinski3} S. Koike, A. Parusi\'nski, 
\emph{ Equivalence relations for two variable real analytic function germs, } arXiv:0801.2650


\bibitem{kuo}  T.-C. Kuo,
   \emph{On classification of  real singularities,}
   Invent. math. \textbf{82} (1985), 257--262.

\bibitem{kurdykaAS} K. Kurdyka,
\emph{Ensembles semi-alg\'ebriques
sym\'etriques par arcs}, Math. Ann. \textbf{ 281} no.3  (1988),  445--462.

\bibitem{aussoisKP} K. Kurdyka, A. Parusi\'nski,
\emph{Arc-symmetric sets and arc-analytic mappings}, 
in 
"Arc Spaces and Additive Invariants in Real Algebraic Geometry",  Panoramas \&
Synth\`eses \textbf{24}, S.M.F. (2007), 33--67.

\bibitem{ens}  C. McCrory, A. Parusi\'nski,
\emph{   Algebraically Constructible Functions}, 
Annales  Sci. \'Ec. Norm. Sup.  
\textbf{30}, 4 (1997), 527--552. 

\bibitem{virtual}  C. McCrory, A. Parusi\'nski,
\emph{Virtual Betti numbers of real algebraic varieties},
Comptes Rendus Acad. Sci. Paris, Ser. I, \textbf{336} (2003),
763--768. 

\bibitem{weight}  C. McCrory, A. Parusi\'nski,
\emph{The weight filtration for real algebraic varieties},
mathArxiv. AG/0807.4203, to appear in  the proceedings of the MSRI workshop "Topology of 
Stratified Spaces".  


\bibitem  {SubFun}   A. Parusi\'nski,
\emph{   Subanalytic functions},
Trans. Amer. Math. Soc.  \textbf{344}, 2 (1994), 583--595. 


\bibitem  {mathannalen}   A. Parusi\'nski,
\emph{Topology of injective endomorphisms of real algebraic sets, }
Math. Annalen, 
\textbf{328} (2004), 353--372.


\end{thebibliography}
\end{document}